\documentclass[11pt,dvips,twoside,letterpaper]{article}
\usepackage{pslatex}
\usepackage{fancyhdr}
\usepackage{graphicx}
\usepackage{geometry}

\def\figurename{Figure} % Replace the colon that normally appears after the Figure number by a period.
\makeatletter
\renewcommand{\fnum@figure}[1]{\figurename~\thefigure.}
\makeatother

\def\tablename{Table} % Replace the colon that normally appears after the Figure number by a period.
\makeatletter
\renewcommand{\fnum@table}[1]{\tablename~\thetable.}
\makeatother

\usepackage{amsmath}
\usepackage{amssymb}
\usepackage{amsfonts}
\usepackage{amsthm,amscd}

\newtheorem{theorem}{Theorem}[section]

\newtheorem{corollary}[theorem]{Corollary}
\newtheorem{proposition}[theorem]{Proposition}
\theoremstyle{definition}
\newtheorem{definition}[theorem]{Definition}

\theoremstyle{remark}

\numberwithin{equation}{section}

\def\P{\mathbb P}

\def\cal{\mathcal}

\begin{document}
%\vskip 0.4in
\title{\bfseries\scshape{LIGHTLIKE OSSERMAN SUBMANIFOLDS OF SEMI-RIEMANNIAN MANIFOLDS}}
\author{\bfseries\scshape Cyriaque Atindogbe$^{a}$\thanks{atincyr@imsp-uac.org}, \bfseries\scshape Oscar Lungiambudila$^{a}$\thanks{lungiambudila@yahoo.fr ; lungiaoscar@imsp-uac.org (Corresponding author)} , \bfseries\scshape Jo\"{e}l Tossa$^{a}$\thanks{joel.tossa@imsp-uac.org ; joel.tossa@uac.imsp.bj }\\
$^{a}$Institut de Math\'ematiques et de Sciences Physiques (IMSP)\\Universit\'e d'Abomey-Calavi\\P.B. 613, Porto-Novo, B\'enin} 
 
\date{}
\maketitle \thispagestyle{empty} \setcounter{page}{1}

% ------- [First Page Running Head] - place it immediately after title! ------
%\thispagestyle{fancy} \fancyhead{}
%\fancyhead[L]{{\LARGE A}frican {\LARGE D}iaspora {\LARGE J}ournal of {\LARGE M}athematics\\
%Volume XX, Number XX, pp. {\thepage--\pageref{lastpage-} (2010)}} % put \label{lastpage-xx} on the last page!
%\fancyhead[R]{ISSN 1539-854X  \\ {\sf www.african-j-math.org}}
%\fancyfoot{}
%\renewcommand{\headrulewidth}{0pt}
%------------------------------------------------------------------------------

\begin{abstract}
In this paper, we study 
Jacobi operators associated to algebraic curvature maps (tensors)
on lightlike submanifolds $M$. We investigate conditions for an induced Riemann curvature tensor to be an algebraic curvature tensor on $M$. We introduce the notion of lightlike Osserman submanifolds and an example of $2$-degenerate Osserman metric is given. Finally we give some results of symmetry properties on lightlike hypersurfaces from Osserman condition. 
\end{abstract}

\noindent {\bf AMS Subject Classification:} 53B25, 53C25, 53C50  

\vspace{.08in} \noindent \textbf{Keywords}: Algebraic curvature tensor, Admissible pair of screens, Pseudo-Jacobi operator, Lightlike Osserman submanifold. 
\section{Introduction}
The curvature tensor is a central concept in differential geometry. 
According to R. Osserman (\cite{Oss}), one could argue that it is a central one. 
But the curvature tensor is in general difficult to deal with and the problem which aims to relate
 algebraic properties of the Riemann curvature tensor to the geometry of the manifold is in general
 difficult to be solved. Many authors study the geometric consequences that followed if various 
natural operators defined in terms of the curvature tensor are assumed to have constant 
eingenvalues on the unit fibre bundle. Osserman has studied the spectral properties of Jacobi 
operator in (\cite{Oss}). This operator has been extensively studied in the Riemannian and the 
pseudo-Riemannian context. E. Garc\'{\i}a-R\'{\i}o, D. N. Kupeli and R. V\'{a}quez-Lorenzo have studied Osserman condition in Pseudo-Riemannian geometry in (\cite{Ga-Ku-Va}). We refer to (\cite{Ga-Ku-Va}) for an extensive bibliography. In the degenerate geometry, C. Atindogbe and K. L. Duggal have stadied Pseudo-Jacobi operators and introduce the Osserman condition on lightlike hypersurfaces in (\cite{Ati-Dug}). In the present paper we extend this study on $r$-degenerate submanifolds and give some characterization results of symmetry properties on lightlike Osserman hypersurfaces. 

Let $(M,g)$ be a semi-Riemannian manifold and $u \in M$. An element $F\in \otimes^{4}T_{u}^{\ast}M$ is said to 
be an {\itshape algebraic curvature map (tensor)} on $T_{u}M$ if it satisfies the following symmetries:
\begin{eqnarray}\label{int1}  
 F(x,y,z,w)=-F(y,x,z,w)=F(z,w,x,y)\crcr
 F(x,y,z,w)+F(y,z,x,w)+F(z,x,y,w)=0
\end{eqnarray}
In semi-Riemannian geometry, the Riemann curvature tensor $R$ is an algebraic curvature tensor on the tangent space $T_{u}M$ for every point $u \in M$.

For an algebraic curvature map (tensor) $F$ on $T_{u}M$, the associated Jacobi operator $J_{F}(x)$ with respect to $x\in T_{u}M$, is the self-adjoint linear map on $T_{u}M$ characterized by identity
\begin{equation}\label{int2} 
 g(J_{F}(x)y,z)=F(y,x,x,z),\quad\forall y,z\in T_{u}M.
\end{equation}
Since $J_{F}(cx)=c^{2}J_{F}(x)$, the natural domains of Jacobi operators 
$J_{F}(\cdotp)$ are the unit pseudo-sphere of unit timelike or unit spacelike vectors 
$$S_{u}^{\pm}(M) := \left\{x \in T_{u}M : g(x,x) = \pm 1\right\}.$$

The tensor $F\in \otimes^{4}T_{u}^{\ast}M$ is said to be a spacelike
(resp. timelike) Osserman tensor on $T_{u}M$ if $Spec\{J_{F}\}$ 
is constant on the pseudo-sphere $S_u^{+}(M)$ (resp $S_u^{-}(M)$ ). 

In degenerate geometry, it is known that the induced metric on an $r$-degenerate submanifold of a semi-Riemannian manifold has a non-trivial kernel, so the relation (\ref{int2}) is not well defined in the usual way. Also, except a totally geodesic case, in general induced Riemann curvature tensors on a lightlike submanifolds are not algebraic curvature tensors. Therefore, in section 3, we study conditions on a lightlike submanifold to have an induced algebraic Riemann curvature tensor.

In section 4, we have extended pseudo-Jacobi operator on
$r$-degenerate submanifolds, by using non-degenerate metric $\tilde{g}$
associated to the degenerate metric $g$ (see Preliminaries). We introduce and study a class of lightlike Osserman submanifolds. Some results are obtained (Theorem \ref{th4.1}, Proposition \ref{th4.2} and Theorem \ref{th4.3}) and We exhibit an example of $2$-degenerate metric verifying Osserman condition.

It is natural to impose condition on semi-Riemannian manifold that its Riemannian curvature
tensor $R$ be parallel, that is, have vanishing covariant differential, $\nabla R$, where   $\nabla$ is the Levi-Civita connection on semi-Riemannian manifold and $R$ is the corresponding curvature tensor. Such a manifold is said to be locally symmetric. This class of manifolds contains one of manifolds of constant curvature. A semi-Riemannian manifold is called semi-symmetric, if $R\cdot R=0$, which is the integrability condition of $\nabla R=0$. The semi-symmetric manifolds have been classified, in Riemannian case, by Szabo in \cite{Sza1} and \cite{Sza2}. A semi-Riemannian manifold is called Ricci semi-symmetric, if $R\cdot Ric = 0$.

In section 5, we are interested to answer to the following question: ``Are conditions $\nabla R = 0$ and $R\cdot R = 0$ equivalent on lightlike hypersurfaces of semi-Riemannian manifolds?'' These equivalences are not true in general. In virtue of result given by Sahin (\cite{Sahin}, Theorem 4.2), we see that the conditions $\nabla R = 0$ and $R\cdot R = 0$ are equivalent on lightlike hypersurfaces of semi-Euclidean spaces under conditions $Ric(\xi,X)=0$ and $A_{N}\xi$ a vector field non-null. In this paper we give an affirmative answer to this question for lightlike Osserman hypersurfaces of semi-Riemannian manifolds of constant sectional curvature, under condition $A_{N}\xi$ a vector field non-null (Corollary \ref{th5.4}). Also, in the same section, we show that lightlike Osserman hypersurfaces of semi-Riemannian manifolds of constant sectional curvature are Ricci semi-symmetric (Theorem \ref{th5.5}). 
%%%%%%%%%%%%%%%%%%%%%%%%%%%%%%%%%%%%%%%%%%%%%%%%%%%%%%%%%%%%%%%%%%%%%%%%%%%%%%%%%%%%%%%%%
\section{Preliminaries}
\subsection{Lightlike submanifolds of semi-Riemannian manifolds}
We follow (\cite{Dug}) for the notations and formulas used in this paper. Let $(\overline{M},\overline{g})$ be an $(m+n)$-dimensional semi-Riemannian manifold of constant
 index $\nu$, $1\leqslant\nu<m+n$ and $M$ be a submanifold of $\overline{M}$ of codimensional $n$. We assume that both $m$ and $n$ are $\geqslant 1$. At a point $u\in M$, we define the orthogonal complement
 $T_{u}M^{\perp}$ of the tangent space $T_{u}M$ by
$$T_{u}M^{\perp}=\{X_{u}\in T_{u}\overline{M} : \overline{g}(X_{u},Y_{u})=0,\ \forall Y_{u}\in T_{u}M\} $$

We put $RadT_{u}M=RadT_{u}M^{\perp}=T_{u}M\cap T_{u}M^{\perp}$. The submanifold $M$ of $\overline{M}$ is 
said to be an r-lightlike submanifold (one suppose that the index of $\overline{M}$ is $\nu\geqslant r$), if the mapping
$$RadTM : u \in M\longrightarrow RadT_{u}M$$
defines a smooth distribution on $M$ of rank $r>0$. We call $RadTM$ the radical distribution on 
$M$. In the sequel, an $r$-lightlike submanifold will simply be called a lightlike submanifold and
 $g$ is lightlike metric, unless we need to specify $r$.\\
Let $S(TM)$ be a screen distribution which is a semi-Riemannian complematary distribution of $Rad(TM)$ in $TM$, that is, 
\begin{equation}\label{pr1}
 TM=RadTM\perp S(TM).
\end{equation}
We consider a screen transversal vector bundle $S(TM^{\perp})$, which is a semi-Riemannian complementary vector bundle of $Rad(TM)$ in $TM^{\perp}$. Since, for any local frame $\{\xi_{i}\}$ of $Rad(TM)$, there exists a local frame $\{N_{i}\}$ of sections with values in the orthogonal complement of $S(TM^{\perp})$ in $S(TM)^{\perp}$ such that $\overline{g}(\xi_{i},N_{j})=\delta_{ij}$ and $\overline{g}(N_{i},N_{j})=0$, it follows that there exists a lightlike transversal vector bundle $ltr(TM)$ locally spanned by $\{N_{i}\}$ (see \cite{Dug}, p144). Let $tr(TM)$ be complementary (but not orthogonal) vector bundle to $TM$ in $T\overline{M}_{|M}$. Then
\begin{equation}\label{pr2}
 tr(TM)=ltr(TM)\perp S(TM^{\perp}),
\end{equation}
\begin{equation}\label{pr3}
T\overline{M}_{|M}=TM\oplus tr(TM)=S(TM)\perp (RadTM\oplus ltr(TM))\perp S(TM^{\perp}).
\end{equation}
Although $S(TM)$ is not unique, it is canonically isomorphic to the factor vector bundle $TM/ RadTM$ (\cite{Kup}).\\
Throughout this paper, we will discuss the dependence (or otherwise) of the results on induced objects and refer to (\cite{Dug}) for their transformation equations. We say that a submanifold $(M,g,S(TM),S(TM^{\perp}))$ of $\overline{M}$ is\\
(1) r-lightlike if $r<min\{m,n\}$;\\
(2) coisotropic if $r=n$, hence, $S(TM^{\perp})=\{0\}$;\\
(3) isotropic if $r=m<n$, hence $S(TM)=\{0\}$;\\
(4) totally lighlike if $r=m=n$, hence $S(TM)=\{0\}=S(TM^{\perp})$.\\
The Gauss and Weingarten equations are
\begin{equation}\label{pr4}
\overline{\nabla}_{X}Y=\nabla_{X}Y+h(X,Y),\quad\forall X,Y\in\Gamma(TM),
\end{equation}
\begin{equation}\label{pr5}
\overline{\nabla}_{X}V=-A_{V}X+\nabla_{X}^{t}V,\quad\forall X\in\Gamma(TM),\ V\in\Gamma(tr(TM)),
\end{equation}
where $\{\nabla_{X}Y,A_{V}X\}$ and $\{h(X,Y),\nabla_{X}^{t}V\}$ belong to $\Gamma(TM)$ and $\Gamma(tr(TM))$, respectively. $\nabla$ and $\nabla^{t}$ are linear connections on $M$ and on the vector bundle $tr(TM)$, respectively. Moreover, we have
\begin{equation}\label{pr6}
\overline{\nabla}_{X}Y=\nabla_{X}Y+h^{l}(X,Y)+h^{s}(X,Y),\quad\forall X,Y\in\Gamma(TM),
\end{equation}
\begin{equation}\label{pr7}
\overline{\nabla}_{X}N=-A_{N}X+\nabla_{X}^{l}N+D^{s}(X,N),\quad\forall X\in\Gamma(TM),\ N\in\Gamma(ltr(TM)),
\end{equation}
\begin{equation}\label{pr8}
\overline{\nabla}_{X}W=-A_{W}X+\nabla_{X}^{s}W+D^{l}(X,W),\quad\forall X\in\Gamma(TM),\ W\in\Gamma(S(TM^{\perp})).
\end{equation}
Denote the projection of $TM$ on $S(TM)$ by $P$. Then, by using (\ref{pr6})-(\ref{pr8}) and taking into account that $\overline{\nabla}$ is a metric connection, we obtain
\begin{equation}\label{pr9}
\overline{g}(h^{s}(X,Y),W)+\overline{g}(Y,D^{l}(X,W))=g(A_{W}X,Y), 
\end{equation}
\begin{equation}\label{pr10}
\overline{g}(D^{s}(X,N),W)=\overline{g}(N,A_{W}X) 
\end{equation}
From the decomposition (\ref{pr1}) of the tangent bundle of lighlike submanifold, we have
\begin{equation}\label{pr11}
\nabla_{X}Y=\stackrel{*}{\nabla}_{X}PY+\stackrel{*}{h}(X,PY),\quad\forall X,Y\in\Gamma(TM),
\end{equation}
\begin{equation}\label{pr12}
\nabla_{X}\xi=-\stackrel{*}{A}_{\xi}X+\stackrel{*}{\nabla^{t}}_{X}\xi,\quad\forall X\in\Gamma(TM),\ \xi\in\Gamma(RadTM).
\end{equation}
By using the abovee equations, we obtain
\begin{equation}\label{pr13}
\overline{g}(h^{l}(X,PY),\xi)=g(\stackrel{*}{A}_{\xi}X,PY),
\end{equation}
\begin{equation}\label{pr14}
\overline{g}(\stackrel{*}{h}(X,PY),N)=g(A_{N}X,PY), 
\end{equation}
\begin{equation}\label{pr15}
\overline{g}(h^{l}(X,\xi),\xi)=0,\qquad\stackrel{*}{A}_{\xi}\xi=0.
\end{equation}
In general, the induced connection $\nabla$ on $M$ is not a metric connection. Since $\overline{\nabla}$ is a metric connection, by using (\ref{pr6}) we get
\begin{equation}\label{pr16}
(\nabla_{X}g)(Y,Z)=\overline{g}(h^{l}(X,Y),Z)+\overline{g}(h^{l}(X,Z),Y),\quad\forall X,Y,Z\in\Gamma(TM).
\end{equation}
However, it is important to note that $\stackrel{*}{\nabla}$ is a metric connection on $S(TM)$.\\
We denote the Riemann curvature tensors of $\overline{M}$ and $M$ by $\overline{R}$ and $R$ respectively. The Gauss equation for $M$ is given by
\begin{eqnarray}\label{pr17}
\overline{R}(X,Y)Z&=&R(X,Y)Z+A_{h^{l}(X,Z)}Y-A_{h^{l}(Y,Z)}X+A_{h^{s}(X,Z)}Y\crcr
&&-A_{h^{s}(Y,Z)}X+(\nabla_{X}h^{l})(Y,Z)-(\nabla_{Y}h^{l})(X,Z)\crcr
&&+D^{l}(X,h^{s}(Y,Z))-D^{l}(Y,h^{s}(X,Z))+(\nabla_{X}h^{s})(Y,Z)\crcr
&&-(\nabla_{Y}h^{s})(X,Z)+D^{s}(X,h^{l}(Y,Z))-D^{s}(Y,h^{l}(X,Z)),
\end{eqnarray}
Therefore,
\begin{eqnarray}\label{pr17.1}
\overline{R}(X,Y,Z,PU)&=&R(X,Y,Z,PU)+\overline{g}(\stackrel{*}{h}(Y,PU),h^{l}(X,Z))-\overline{g}(\stackrel{*}{h}(X,PU),h^{l}(Y,Z)) \crcr 
&&+\overline{g}(h^{s}(Y,PU),h^{s}(X,Z))-\overline{g}(h^{s}(X,PU),h^{s}(Y,Z))
\end{eqnarray}
for any $X,Y,Z\in\Gamma(TM)$. Note that for the coisotropic, isotropic and totally lighlike submanifolds, in (\ref{pr17}), we have $h^{s}=0$, $h^{l}=0$ and $h^{l}=h^{s}=0$, respectively. 
\subsection{Pseudo-inversion of $r$-degenerate metrics}
In this section we mention (case 1 or 2) to refer to $r$-lightlike submanifolds with $0<r<min\{m,n\}$ or coisotropic submanifolds and also (case 3 or 4) for isotropic submanifolds or totally lightlike submanifolds. We recall from (\cite{Lun}) the following result. Consider on $M$ the local frames $\{\xi_{i}\}$ and $\{N_{i}\}$ of sections of $RadTM$ and $ltr(TM)$ satisfying $\overline{g}(N_{i},\xi_{j})=\delta_{ij}$. Consider on $M$ the 1-forms $\eta_{i},i=1,...,r$ defined by $\eta_{i}(\cdotp)=\overline{g}(N_{i},\cdotp)$. Any vector field $X$ on $M$ is expressed on a coordinate neighbourhood $\mathcal{U}$ as follows,
\begin{equation}\label{pr18} 
X=PX+\sum_{i=1}^{r}\eta_{i}(X)\xi_{i}\qquad (case\ 1~or~2)
\end{equation}
\begin{equation}\label{pr19} 
X=\sum_{i=1}^{m}\eta_{i}(X)\xi_{i}\qquad\qquad (case\ 3~or~4) 
\end{equation} 
Now, we define $\flat_{g}$ by
\begin{eqnarray*}
 \flat_{g}:\Gamma(TM)&\longrightarrow& \Gamma(T^{*}M) \crcr
           X & \longmapsto & X^{\flat_{g}}
\end{eqnarray*}
such that for all $Y\in\Gamma(TM)$,
\begin{eqnarray} 
&&X^{\flat_{g}}(Y)=g(X,Y)+\sum_{i=1}^{r}\eta_{i}(X)\eta_{i}(Y)\qquad (case\ 1~or~2)\label{pr20} \\
&&X^{\flat_{g}}(Y)=\sum_{i=1}^{m}\eta_{i}(X)\eta_{i}(Y)\qquad\qquad (case\ 3~or~4)\label{pr21}
\end{eqnarray}
The map $\flat_{g}$ is an isomorphism of $\Gamma(TM)$ onto $\Gamma(T^{*}M)$, its inverse is denoted $\sharp_{g}$.  For $X\in\Gamma(TM)$ (resp. $\omega\in\Gamma(T^{*}M)$), $X^{\flat_{g}}$ (resp. $\omega^{\sharp_{g}}$) is called the dual 1-form of $X$ (resp. the dual vector field of $\omega$) with respect to the degenerate metric $g$.\\
We define a $(0,2)$-tensor $\tilde{g}$ by, for any $X,Y\in\Gamma(TM)$,
\begin{equation}\label{pr22}
\tilde{g}(X,Y)= X^{\flat_{g}}(Y)= g(X,Y)+\sum_{i=1}^{r}\eta_{i}(X)\eta_{i}(Y)\qquad (case\ 1~or~2) 
\end{equation}
and
\begin{equation}\label{pr23}
\tilde{g}(X,Y)=X^{\flat_{g}}(Y)= \sum_{i=1}^{m}\eta_{i}(X)\eta_{i}(Y)\qquad (case\ 3~or~4) 
\end{equation}
Clearly, $\tilde{g}$ defines a non-degenerate metric on $M$. Also, observe that $\tilde{g}$ coincides with $g$ if the latter is non-degenerate. The $(0,2)$-tensor $\tilde{g}^{-1}$, inverse of $\tilde{g}$ is called the pseudo-inverse of $g$. Let consider the local quasi-orthonormal field of frames $\{\xi_{1},...,\xi_{r},X_{r+1},...,X_{m}\}$ and $\{\xi_{1},...,\xi_{r}\}$ on lightlike submanifold $M$ with respect to the decompositions $TM=S(TM)\perp RadTM$ (case 1 or 2 ) and $TM=RadTM$ (case 3 or 4). Using relations (\ref{pr22}) and (\ref{pr23}), we have\\
$$\tilde{g}(\xi_{i},\xi_{j})=\delta_{ij},~1\leqslant i,j\leqslant r\quad and\quad   \tilde{g}(X_{i},X_{j})=g_{ij},~r+1\leqslant i,j\leqslant m,\quad (Case 1~ or~ 2).$$
$$\tilde{g}(\xi_{i},\xi_{j})=\delta_{ij},\quad 1\leqslant i,j\leqslant m,\quad (Case 3~ or~ 4).$$
%%%%%%%%%%%%%%%%%%%%%%%%%%%%%%%%%%%%%%%%%%%%%%%%%%%%%%%%%%%%%%%%%%%%%%%%%%%%%%%%%%%%
\section{Algebraic Riemann cuvature tensors}
Contrary to non-degenerate hypersurfaces, the induced Riemann curvature on lightlike submanifold $(M, g, S(T M),S(TM^{\perp}))$ may not have an algebraic curvature tensor. For this, we have the following results.
\begin{proposition}\label{th3.1}
Let $M$ be an isotropic submanifold or a totally lightlike submanifold of a semi-Riemannian manifold $\overline{M}$. Then, the induced Riemann curvature $R$ on $M$ is an algebraic curvature tensor.
\end{proposition}
\textbf{Proof:} Since $TM=RadTM$, we have, $R(X,Y,Z,U)=0,~~\forall X,Y,Z,U\in\Gamma(TM)$.$\qquad\square$
%%%%%%%%%%%%%%%%%%%%%%
\begin{theorem}\label{th3.2}
Let $M$ be an $r$-lightlike submanifold with $r<min\{m,n\}$ or a coisotropic submanifold of a semi-Riemannian manifold $\overline{M}$, such that the radical distribution $RadTM$ is integrable. If the induced Riemann curvature tensor of $M$ is an algebraic curvature tensor, then at least one of the following holds\\
(1) $h^{l}(X,Y)=0$, for any $X,Y\in\Gamma(TM)$.\\
(2) $A_{N}\xi\in\Gamma(RadTM)$, for any $\xi\in\Gamma(RadTM),~N\in\Gamma(ltr(TM))$\\
(3) $\overline{\nabla}_{\xi}N\in\Gamma(tr(TM))$, for any $\xi\in\Gamma(RadTM),~N\in\Gamma(ltr(TM))$.
\end{theorem}
\textbf{Proof:} Suppose that the radical distribution $RadTM$ is integrable. In virtue of Theorem 2.7 of \cite{Dug}, p162, we have 
\begin{equation}\label{eq3.1}
 h^{l}(PX,\xi),\quad\forall\xi\in\Gamma(RadTM),~X\in\Gamma(TM)
\end{equation}
Moreover, suppose that the induced curvature $R$ of $M$ is algebraic. Since $R(X,\xi,Y,Z)=-R(X,\xi,Z,Y)$, for any $\xi\in\Gamma(RadTM)$ and $X,Y,Z\in\Gamma(S(TM))$, from relations (\ref{pr17.1}) and (\ref{eq3.1}), we have
\begin{equation}\label{eq3.2}
\overline{R}(X,\xi,Y,Z)-\overline{g}(\stackrel{*}{h}(\xi,Z),h^{l}(X,Y))=
-\overline{R}(X,\xi,Z,Y)+\overline{g}(\stackrel{*}{h}(\xi,Y),h^{l}(X,Z)).
\end{equation}
Using symmetry of $h^{l}$ and since the Riemann curvature $\overline{R}$ of $\overline{M}$ is algebraic curvature tensor, the relation (\ref{eq3.2}) leads to
\begin{eqnarray*}
\overline{g}(\stackrel{*}{h}(\xi,Z),h^{l}(X,Y))&=&-\overline{g}(\stackrel{*}{h}(\xi,Y),h^{l}(X,Z))=-\overline{g}(\stackrel{*}{h}(\xi,Y),h^{l}(Z,X))\\
&=&\overline{g}(\stackrel{*}{h}(\xi,X),h^{l}(Z,Y))=\overline{g}(\stackrel{*}{h}(\xi,X),h^{l}(Y,Z))\\
&=&-\overline{g}(\stackrel{*}{h}(\xi,Z),h^{l}(Y,X))=-\overline{g}(\stackrel{*}{h}(\xi,Z),h^{l}(X,Y)).
\end{eqnarray*}
Hence,
\begin{equation}\label{eq3.3}
\overline{g}(\stackrel{*}{h}(\xi,Z),h^{l}(X,Y))=0,\quad\forall\xi\in\Gamma(RadTM),~X,Y,Z\in\Gamma(S(TM)).
\end{equation}
Since $\overline{g}$ is a non-degenerate metric, by using (\ref{eq3.1}), we infer from (\ref{eq3.3}) that $\stackrel{*}{h}(\xi,PZ)=0$ or $h^{l}(X,Y)=0$, for any $X,Y,Z\in\Gamma(TM)$. Now, assume that in (\ref{eq3.3}) there exist $X_{0},Y_{0}\in\Gamma(S(TM))$ such that $h^{l}(X_{0},Y_{0})\neq0$. Then $\stackrel{*}{h}(\xi,PZ)=0,~\forall \xi\in\Gamma(RadTM),~Z\in\Gamma(TM)$. This leads to the following
\begin{equation}\label{eq3.4}
\overline{g}(\stackrel{*}{h}(\xi,PZ),N)=g(A_{N}\xi,PZ)=0,~\forall N\in\Gamma(ltr(TM)),~Z\in\Gamma(TM),
\end{equation}
that is $A_{N}\xi\in\Gamma(RadTM),~\forall \xi\in\Gamma(RadTM),~N\in\Gamma(ltr(TM))$. Moreover, if $A_{N}\xi=0$, we have  
\begin{equation}\label{eq3.5}
\overline{g}(\stackrel{*}{h}(\xi,PZ),N)=\overline{g}(\overline{\nabla}_{\xi}PZ,N)=
-\overline{g}(\overline{\nabla}_{\xi}N,PZ)=0.
\end{equation}
This lead to $\overline{\nabla}_{\xi}N\in\Gamma(tr(TM))$, for any $\xi\in\Gamma(RadTM),~N\in\Gamma(ltr(TM))$.$\qquad\square$
%%%%%%%%%%%%%%%%%%%
\begin{corollary}\label{th3.3}
Let $M$ be a coisotropic submanifold of $\overline{M}$, such that the radical distribution $RadTM$ is integrable. If the induced Riemann curvature tensor of $M$ is an algebraic curvature tensor, then at least one of the following holds\\
(1) $M$ is totally geodesic.\\
(2) $A_{N}\xi\in\Gamma(RadTM)$, for any $\xi\in\Gamma(RadTM),~N\in\Gamma(ltr(TM))$\\
(3) $\overline{\nabla}_{\xi}N\in\Gamma(ltr(TM))$, for any $\xi\in\Gamma(RadTM),~N\in\Gamma(ltr(TM))$.
\end{corollary}
If $(M,g)$ is a lightlike hypersurface, the radical distribution $RadTM=TM^{\perp}$ and it is integrable. Also the shape operator $A_{N}$ is $\Gamma(S(TM))$-valued. So, the following holds.
\begin{corollary}\label{th3.4}
Let $(M,g)$ be a lightlike hypersurface of a semi-Riemannian manifold $(\overline{M},\overline{g})$. If the induced Riemann curvature tensor of $M$ is an algebraic curvature tensor, then at least one of the following holds\\
(1) $M$ is totally geodesic.\\
(2) $\overline{\nabla}_{\xi}N\in\Gamma(tr(TM))$, for any $\xi\in\Gamma(RadTM),~N\in\Gamma(tr(TM))$.
\end{corollary}
Let's consider a family of coisotropic submanifolds $M$ such that the local second fundamental forms of the screen distribution $S(TM)$ are related with the local second fundamental forms of $M$ as follows:
\begin{equation}\label{eq3.6}
\stackrel{*}{h}_{i}(X,PY)=\varphi_{i}h_{i}^{l}(X,PY),\quad\forall X,Y\in\Gamma(TM_{|\cal{U}}),
\end{equation}
where each $\varphi_{i}$ is a conformal smooth function on a coordinate neighbourhood $\cal{U}$ in $M$. The following result holds.
\begin{theorem}\label{th3.5} Let $(M,g,S(TM)$ be a coisotropic submanifold of a semi-Riemannian manifold $(\overline{M},\overline{g})$, such that (\ref{eq3.6}) holds and the radical distribution $RadTM$ is integrable. Then, the induced Riemann curvature $R$ of $M$ defines an algebraic curvature tensor if the following holonomy condition is satisfied
$$\overline{R}(X,PY)RadTM\subset RadTM,\quad\forall X,Y\in\Gamma(TM).$$
\end{theorem}
\textbf{Proof:} Consider $M$ such that $RadTM$ is integrable. Since $\stackrel{*}{h}_{i}(X,PY)=\varphi_{i}h_{i}^{l}(X,PY)$, from (\ref{pr17.1}) we have , for any $X,Y,Z,W\in \Gamma(TM)$,
\begin{eqnarray*}
 R(X,Y,Z,PU)&=&\overline{R}(X,Y,Z,PU)+\sum_{i=1}^{r}\varphi_{i}\big\{h_{i}^{l}(X,PU)h_{i}^{l}(Y,Z)-h_{i}^{l}(Y,PU)h_{i}^{l}(X,Z)\big\}\\
            &=& \overline{R}(X,Y,Z,PU)+\mathcal A(X,Y,Z,PU),
\end{eqnarray*}
where $$\mathcal A(X,Y,Z,U)=\sum_{i=1}^{r}\varphi_{i}\big\{h_{i}^{l}(X,U)h_{i}^{l}(Y,Z)-h_{i}^{l}(Y,U)h_{i}^{l}(X,Z)\big\},~\forall X,Y,Z,U\in \Gamma(TM).$$
It is straightforward that $\mathcal A$ verifies the algebraic symmetries of (\ref{int1}). So, $R(X,Y,Z,PU)$ has the required symmetries. On the other hand, for any $X,Y,Z\in \Gamma(TM)$ and $\xi\in RadTM$, we have $R(X,Y,Z,\xi)=-R(Y,X,Z,\xi)=0$. Also, since $h_{i}^{l}(PX,\xi)=0,~\forall X\in\Gamma(TM),\xi\in\Gamma(RadTM)$, we have
$R(Z,\xi,X,Y)=R(Z,\xi,X,PY)=\overline{R}(Z,\xi,X,PY)=-\overline{R}(X,PY,\xi,Z)=0$. This completes the proof.$\qquad\square$\vspace{0.2cm}\\
It is easy to see that a lightlike hypersurface verifying (\ref{eq3.6}) is said to be locally (globally) screen conformal (see \cite{Ati-Dug2}). So, the following holds.
\begin{corollary}\label{th3.6}
Let $(M,g,S(TM))$ be a locally screen conformal lightlike hypersurface of a semi-Riemannian manifold $(\overline{M},\overline{g})$ with ambient holonomy condition
$$\overline{R}(X,PY)RadTM\subset RadTM,\quad\forall X,Y\in\Gamma(TM).$$
Then, the induced Riemann curvature $R$ of $M$ defines an algebraic curvature tensor.
\end{corollary}

%%%%%%%%%%%%%%%%%%%%%%%%%%%%%%%%%%%%%%%%%%%%%%%%%%%%%%%%%%%%%%%%%%%%%%%%%%%%%%
\section{Pseudo-Jacobi operators and lightlike Osserman submanifolds}
Let's start by intrinsic interpretation of relation (\ref{int2}) which in pseudo-Riemannian setting characterizes the Jacobi operator $J_{R}(\cdotp)$ associated to an algebraic curvature map (tensor) $R\in \otimes^{4}T_{u}^{\ast}M$, $u\in M$. Indeed, for $x\in S_{u}^{+}(M)$ (or $x\in S_{u}^{-}(M)$), $y$, $w$ in $T_{u}M$, we have,
\begin{equation}\label{eq4.11} 
 J_{R}(x)y=R(y,x,x,\bullet)^{\sharp}
\end{equation}
that is
\begin{equation}\label{eq4.12} 
 (J_{R}(x)y)^{\flat}(w)=R(y,x,x,w)
\end{equation} 
where $\flat$ and $\sharp$ are the usual natural isomorphisms between $T_{u}M$ and its dual $T_{u}^{*}M$, with respect to non-degenerate metric $g$. For degenerate setting, let's consider the associate non-degenerate metric $\tilde{g}$ of $g$ defined by relations (\ref{pr22}) and (\ref{pr23}), and denote by $\flat_{g}$ and $\sharp_{g}$ the natural isomorphisms with respect to the metric $\tilde{g}$. Thus, equivalently the above relations can be written in the form:
\begin{equation}\label{eq4.13} 
\tilde{g}(J_{R}(x)y,w)=R(y,x,x,w)
\end{equation}
in which $J_{R}(x)$ is well defined. This leads to the following definition.
\begin{definition}[Pseudo-Jacobi Operator]
Let $(M,g,S(TM),S(TM^{\perp}))$ be a lightlike submanifold of a semi-Riemannian manifold $(\overline{M},\overline{g})$, $u\in M$, $x\in S_{u}^{+}(M)$ (or $x\in S_{u}^{-}(M)$) and $R\in \otimes^{4}T_{u}^{\ast}M$ an algebraic curvature map (tensor) on $T_{u}M$. By {\itshape pseudo-Jacobi operator} associated to $R$ with respect to $x$, we call the self-adjoint linear map $J_{R}(x)$ on $T_{u}M$ defined by
\begin{equation}\label{eq4.14} 
 J_{R}(x)y=R(y,x,x,\bullet)^{\sharp_{g}}
\end{equation}
or equivalently
\begin{equation}\label{eq4.15} 
 (J_{R}(x)y)^{\flat_{g}}(w)=R(y,x,x,w)
\end{equation}
where $\flat_{g}$ and $\sharp_{g}$ denote the natural isomorphisms between $T_{u}M$ and its dual $T_{u}^{*}M$, with respect to non-degenerate metric $\tilde{g}$ of $g$. 
\end{definition}
\textbf{Remark:} For the cases of isotropic submanifold and totally lightlike submanifold, the induced Riemannian curvature tensor $R$ vanish identically on $T_{u}M$ for any $u\in M$. Therefore, the pseudo-Jacobi operator $J_{R}(\cdot)$ associated to $R$ vanish identically on $S_{u}^{\pm}(M)$.\vspace{0.2cm}

It is known by approach developed in (\cite{Dug}) that, the extrinsic geometry of lightlike submanifolds depends on a choice of screen distribution $S(TM)$ and screen transversal vector bundle $S(TM^{\perp})$. Since these screens are not uniquely determined, a well defined concept of Jacobi condition is not possible for an arbitrary lightlike submanifold of a semi-Riemannian manifold. Thus, one must look for a class of pair $\{S(TM),S(TM^{\perp})\}$ of screens for which the induced Riemann curvature and associated pseudo-Jacobi operator have the desired symmetries and properties. In short, we precise the following. 
\begin{definition}
A pair $\{S(TM),S(TM^{\perp})\}$ of screen distribution and screen transversal vector bundle on lightlike submanifold $M$ of a semi-Riemannian manifold $\overline{M}$ is said to be \textit{admissible} if its associated induced Riemann curvature $R$ is an algebraic curvature tensor.
\end{definition}
Note that for case of coisotropic submanifold, since $S(TM^{\perp})=\{0\}$, we'll use the concept \textit{admissible screen distribution} $S(TM)$.\\
\textbf{Example }\\
(1) It is obvious that on totally geodesic lightlike submanifold, any pair of screens is admissible.\\
(2) In virtue of Theorem \ref{th3.5}, any coisotropic submanifold $(M,g,S(TM))$ of semi-Riemannian manifold of constant sectional curvature $\overline{M}(c)$ with integrable radical distribution, satisfying relation (\ref{eq3.6}), admits an admissible screen distribution.
\begin{definition}
A lightlike submanifold $(M,g)$ of a semi-Riemannian manifold $(\overline{M},\overline{g})$ of constant index is called timelike (resp. spacelike) Osserman at $u\in M$ if for each admissible pair of screens  $\{S(TM),S(M^{\perp})\}$ and associate induced Riemann curvature $R$, the characteristic polynomial of $J_{R}(x)$ is independent of $x\in S_{u}^{-}(M)$ (resp. $x\in S_{u}^{+}(M)$). Moreover, if this holds at each $u\in M$, then $(M,g)$ is called pointwise Osserman (or Osserman). If this holds independently of the point $u\in M$, $(M,g)$ is called globally Osserman.
\end{definition}
Note that the above definition of Osserman condition extends the definition given in (\cite{Ati-Dug}) for lightlike hypersurface and is independent on the choice of admissible pair of screens. Also, we can show that a lightlike submanifold $(M,g)$ being timelike Osserman at $p\in M$ is equivalent to $(M,g)$ being spacelike Osserman at $p$.\vspace{0.2cm}\\
According to the above remark, we have the following result.
\begin{theorem}
 Let $(M,g)$ be an isotropic submanifold or a totally lightlike submanifold of a semi-Riemannian manifold $(\overline{M},\overline{g})$. Then $(M,g)$ is globally Osserman.
\end{theorem}
In the following example, we exhibit a $2$-degenerate metric which satisfies Osserman condition. We have considered this metric as an induced metric on a coisotropic submanifold.\newpage
\hspace{-0.6cm}\textbf{Basic example}\\
Let $(x,y)=(x_{0},...,x_{p},y_{0},...,y_{p})$ be the usual coordinates on $\mathbb{R}^{2p+2}$. Let $f=f(x_{1},...,x_{p})$ and $h=h(x_{1},...,x_{p})$ be the smooth functions on an open subset $\mathcal{O}\subset\mathbb{R}^{p}$. We define with respect to the natural field of frames $\{\frac{\partial}{\partial x_{0}},...,\frac{\partial}{\partial x_{p}},\frac{\partial}{\partial y_{0}},...,\frac{\partial}{\partial y_{p}}\}$ a $2$-degenerate metric $g_{(f,h)}$ on $M=\mathbb{R}\times\mathcal{O}\times\mathbb{R}^{p+1}$ by 
\begin{equation}\label{eqex1.1}
g_{(f,h)}=\sum_{i=1}^{p}(\frac{\partial f}{\partial x_{i}}dx_{0}dx_{i}+\frac{\partial h}{\partial x_{i}}dx_{i}dy_{0})+\sum_{i,j=1}^{p}\{(\frac{\partial f}{\partial x_{i}}\frac{\partial f}{\partial x_{j}}+\frac{\partial h}{\partial x_{i}}\frac{\partial h}{\partial x_{j}})dx_{i}dx_{j}+\delta_{ij}dx_{i}dy_{j}\}.
\end{equation}
The $2$-degenerate manifold $(M,g_{(f,h)})$ arise as a lightlike submanifold in a  $(p+4)$-dimensional semi-Euclidean space $\overline{M}$. Let $\{u_{0},...,u_{p},v_{0},...,v_{p},w_{1},w_{2}\}$ be a basis for a space $\overline{M}$. Define an semi-Euclidean metric $\overline{g}$ of signature $(p+2,p+2)$ on $\overline{M}$ by setting
$$\overline{g}(u_{i},u_{j})=0= \overline{g}(v_{i},v_{j}),\quad
\overline{g}(u_{0},v_{j})=0=\overline{g}(u_{i},v_{0}),\quad0\leqslant i,j\leqslant p.$$
$$\overline{g}(u_{0},w_{1})=1\quad\overline{g}(u_{0},w_{2})=0,\quad
\overline{g}(u_{i},w_{1})=0=\overline{g}(u_{i},w_{2}),\quad 1\leqslant i\leqslant p.$$
$$\overline{g}(v_{0},w_{1})=0\quad\overline{g}(v_{0},w_{2})=1,\quad
\overline{g}(v_{i},w_{1})=0=\overline{g}(v_{i},w_{2}),\quad 1\leqslant i\leqslant p.$$
$$\overline{g}(u_{i},v_{j})=\delta_{ij},\quad \overline{g}(w_{i},w_{j})=\delta_{ij},\quad
 \quad1\leqslant i,j\leqslant p.$$
Let consider the application
\begin{equation}\label{eqex1.2}
F(x,y)=x_{0}~u_{0}+...+x_{p}~u_{p}+y_{0}~v_{0}+...+y_{p}~v_{p}+f~w_{1}+h~w_{2}.
\end{equation}
$F(x,y)$ defines an embedding of $M$ in $\overline{M}$ and $g_{(f,h)}$ is the induced metric on the embedded submaifold $M$.\vspace{0.2cm}\\
FACT 1. By direct calculation using (\ref{eqex1.2}), the tangent space $TM$ is defined by
\begin{eqnarray}
TM&=&Span\{\partial_{0}^{x}=u_{0},~\partial_{1}^{x}=u_{1}+\partial_{1}^{x}f~w_{1}+\partial_{1}^{x}h~w_{2},...,\crcr
&&\qquad\quad\partial_{p}^{x}=u_{p}+\partial_{p}^{x}f~w_{1}+\partial_{p}^{x}h~w_{2},
\partial_{0}^{y}=v_{0},~\partial_{1}^{y}=v_{1},...,\partial_{p}^{y}=v_{p}\},
\end{eqnarray}
where $\partial_{i}^{x}=\frac{\partial}{\partial x_{i}}$ and $\partial_{i}^{y}=\frac{\partial}{\partial y_{i}}$.\\
The radical distribution $RadTM$ of rank $2$ is given by 
\begin{equation}\label{eqex1.3}
RadTM=Span
\big\{\xi_{1}=\partial_{0}^{x}-\sum_{i=1}^{p}\partial_{i}^{x}f~\partial_{i}^{y}~,
~\xi_{2}=\partial_{0}^{y}-\sum_{i=1}^{p}\partial_{i}^{x}h~\partial_{i}^{y}\big\}.
\end{equation}
$M$ is a coisotropic submanifold of a semi-Euclidean space $\overline{M}$. The lightlike transversal vector bundle $ltr(TM)$ of $M$ is given by
\begin{equation}\label{eqex1.4}
ltr(TM)=Span
\big\{N_{1}=w_{1}-\frac{1}{2}\xi_{1}~,~N_{2}=w_{2}-\frac{1}{2}\xi_{2}\big\}.
\end{equation}
The corresponding screen distribution $S(TM)$ for the above $ltr(TM)$ is given by
\begin{equation}\label{eqex1.5}
S(TM)=\big\{U_{1},...,U_{p},V_{1},...,V_{p}\big\},
\end{equation}
where $U_{i}=\partial_{i}^{x}-\partial_{i}^{x}f~\partial_{0}^{x}-\partial_{i}^{x}h~\partial_{0}^{y}$ and $V_{i}=\partial_{i}^{y}$.\vspace{0.2cm}\\
FACT 2.  Let's consider on $\overline{M}$ a local quasi-orthogonal frame $\{\xi_{1},\xi_{2},U_{i},V_{i},N_{1},N_{2}\}_{1\leqslant i\leqslant p}$ such that $\{\xi_{1},\xi_{2},U_{i},V_{i}\}_{1\leqslant i\leqslant p}$ is a local frame on $M$ with respect to the decomposition (\ref{pr1}). Using the metric $\overline{g}$, we have
$$\overline{\nabla}_{\xi_{1}}^{\xi_{2}}=\overline{\nabla}_{\xi_{2}}^{\xi_{1}}=0,\quad\overline{\nabla}_{V_{i}}^{\xi_{1}}=\overline{\nabla}_{V_{i}}^{\xi_{2}}=0,~1\leqslant i\leqslant p$$
$$\overline{\nabla}_{U_{i}}^{\xi_{1}}=-\sum_{j=1}^{p}\partial_{i}^{x}\partial_{j}^{x}f~V_{j},\quad
\overline{\nabla}_{U_{i}}^{\xi_{2}}=-\sum_{j=1}^{p}\partial_{i}^{x}\partial_{j}^{x}h~V_{j},~1\leqslant i\leqslant P.$$
Thus, since $h_{i}^{l}(X,Y)=-\overline{g}(\overline{\nabla}_{X}^{\xi_{i}},Y)$, we obtain that the non-vanishing components of $h^{l}$ are
\begin{equation}\label{eqex1.6}
h_{1}^{l}(U_{i},U_{j})=\partial_{i}^{x}\partial_{j}^{x}f,\quad h_{2}^{l}(U_{i},U_{j})=\partial_{i}^{x}\partial_{j}^{x}h,~1\leqslant i,j\leqslant p.
\end{equation}
Also, by straightforward calculation, using the Gauss equation, we obtain that the only non-vanishing components of induced connection $\nabla$ on $M$ are
\begin{eqnarray}\label{eqex1.7}
&& \nabla_{U_{i}}^{U_{j}}=-\frac{1}{2}\partial_{i}^{x}\partial_{j}^{x}f~\xi_{1}
-\frac{1}{2}\partial_{i}^{x}\partial_{j}^{x}h~\xi_{2}
-\sum_{k=1}^{p}\big(\partial_{i}^{x}\partial_{j}^{x}f~\partial_{k}^{x}f
+\partial_{i}^{x}\partial_{j}^{x}h~\partial_{k}^{x}h\big)V_{k},\crcr
&& \nabla_{U_{i}}^{\xi_{1}}=-\sum_{j=1}^{p}\partial_{i}^{x}\partial_{j}^{x}f~V_{j},\quad
\nabla_{U_{i}}^{\xi_{2}}=-\sum_{j=1}^{p}\partial_{i}^{x}\partial_{j}^{x}h~V_{j},\quad1\leqslant i,j\leqslant p.
\end{eqnarray}
FACT 3. By direct calculation, the only non-vanishing components of the induced Riemannian curvature tensor on $M$ are given by
\begin{equation}\label{eqex1.8}
R(U_{i},U_{j})U_{k}=\frac{1}{2}\sum_{l=1}^{p}\big\{ f_{;ik}f_{;jl}-f_{;jk} f_{;il}
+h_{;ik}h_{;jl}-h_{;jk}h_{;il}\big\}V_{l}.
\end{equation}
\begin{equation}\label{eqex1.8}
R(U_{i},U_{j},U_{k},U_{l})=\frac{1}{2}\big\{ f_{;ik}f_{;jl}-f_{;jk} f_{;il}
+h_{;ik}h_{;jl}-h_{;jk}h_{;il}\big\} ,
\end{equation}
where $f_{;ij}=\partial_{i}^{x}\partial_{j}^{x}f$ and $h_{;ij}=\partial_{i}^{x}\partial_{j}^{x}h$.\vspace{0.2cm}\\
FACT 4. In virtue of FACT 3, we can see that for any $X,Y,Z\in\Gamma(TM)$, $R(X,Y)Z\in\Gamma(S(TM))$. Thus for any $X\in S^{\pm}(M)$, $Y,W\in\Gamma(TM)$, we have
$$\tilde{g}(J_{R}(X)Y,W)=R(Y,X,X,W)=g(R(Y,X)X,W)=\tilde{g}(R(Y,X)X,W).$$
Since the metric $\tilde{g}$ is non-degenerate, we infer that, for any $X\in S^{\pm}(M)$
\begin{equation}\label{eqex1.9}
J_{R}(X)Y=R(Y,X)X.
\end{equation}
Let $X=X_{a}\xi_{1}+X_{b}\xi_{2}+X_{1}U_{1}+...+X_{p}U_{p}+X_{p+1}V_{1}+...+X_{2p}V_{p}$ the tangent vector field on $M$, by straightforward calculation, with respect to the local quasi-orthogonal frame $\{\xi_{1},\xi_{2},U_{i},V_{i}\}_{1\leqslant i\leqslant p}$ on $M$, we have
$$J_{R}(X)\xi_{1}=0,\quad J_{R}(X)\xi_{2}=0,\quad J_{R}(X)V_{i}=0,$$
\begin{eqnarray*}
J_{R}(X)U_{i}&=&\sum_{j,k=1}^{p}X_{j}X_{k}R(U_{i},U_{j})U_{k}\\
 &=&\sum_{l=1}^{p}\big\{\frac{1}{2}\sum_{j,k=1}^{p}X_{j}X_{k}\big(
f_{;ik}f_{;jl}-f_{;jk} f_{;il}+h_{;ik}h_{;jl}-h_{;jk}h_{;il}\big)\big\}V_{l}\\
 &=& \sum_{l=1}^{p}\Phi_{li}V_{l},\quad 1\leqslant i \leqslant p,
\end{eqnarray*}
where 
$$\Phi_{li}=\frac{1}{2}\sum_{j,k=1}^{p}X_{j}X_{k}\big(
f_{;ik}f_{;jl}-f_{;jk} f_{;il}+h_{;ik}h_{;jl}-h_{;jk}h_{;il}\big),\quad 1\leqslant l,i \leqslant p.$$
Thus, the pseudo-Jacobi operator is given by, for any $X\in S^{\pm}(M)$,
\begin{equation}\label{eqex1.10}
 J_{R}(X)=\left( 
\begin{array}{ccccccc}
& & & & & &\\
&O_{p+2,2} & \vdots & O_{p+2,p} & \vdots & O_{p+2,p} & \\ 
% & \vdots &  \\
&\cdots & \cdots & \cdots& \cdots&\cdots&\\
%& \vdots & \\
&O_{p,2} & \vdots & A_{p,p} & \vdots & O_{p,p} & \\ 
& & & & & &\\
 \end{array}\right)
\end{equation}
where the submatrix $A_{p,p}=\big(\Phi_{li}\big)_{1\leqslant l,i\leqslant p}$. It follows from the expression of $J_{R}(X)$ that its characteristic polynomial satisfies
\begin{equation}\label{eqex1.11}
P_{\lambda}(J_{R}(X))=det(J_{R}(X)-\lambda I_{2p+2})=\lambda^{2p+2}.
\end{equation}
Thus, all eigenvalues are zero. This proves that the $2$-degenerate submanifold $(M,g_{(f,h)})$ is globally Osserman.
\begin{definition}\cite{Dug-Jin}
 A lightlike submanifold $(M,g)$ of a semi-Riemannian manifold
$(\overline{M},\overline{g})$ is said to be totally umbilical in $\overline{M}$ if there is a smooth transversal vector field $\mathbf H\in \Gamma(tr(TM))$ on $M$, called the \textit{transversal curvature vector field} of $M$, such that for any $X,Y\in\Gamma(TM)$,
\begin{equation}\label{eq4.1}
 h(X,Y)=\mathbf Hg(X,Y).
\end{equation}
\end{definition}
Using (\ref{pr6}), (\ref{pr9}) and (\ref{eq4.1}), it is easy to see that $M$ is \textit{totally umbilical} if and only if on each coordinate neighborhood $\mathcal{U}$, there exist smooth vector fields $H^{l}\in\Gamma(ltr(TM))$ and
$H^{s}\in\Gamma(S(TM^{\perp}))$ such that
\begin{eqnarray}\label{eq4.2}
 &h^{l}(X,Y)=H^{l}g(X,Y),\quad D^{l}(X,W)=0&\crcr
&h^{s}(X,Y)=H^{s}g(X,Y),\quad\forall X,Y\in\Gamma(TM),~~W\in\Gamma(S(TM^{\perp})).&
\end{eqnarray}
Note that for the case of lightlike hypersurface, the relations (\ref{eq4.1}) and (\ref{eq4.2}) are equivalent to
\begin{equation}\label{eq4.2.1}
 B(X,Y)=\rho g(X,Y),\quad\forall X,Y\in\Gamma(TM),
\end{equation}
where $\rho$ is the smooth function on $\mathcal{U}\subset M$ and $B$, the local second fundamental form of $M$.
\begin{definition}
Let $(M,g)$ be an $r$-lightlike submanifold with $r<min\{m,n\}$ or a coisotropic submanifold of a $(m+n)$-dimensional semi-Riemannian manifold $(\overline{M},\overline{g})$. We say that the screen distribution $S(TM)$ is \textit{totally umbilical} if for any section $N$ of $ltr(TM)$ on a coordinate neighbourhood $\mathcal{U}\subset M$, there exists a smooth function $\lambda$ on $\mathcal{U}$ such that 
\begin{equation}\label{eq4.3}
\overline{g}(\stackrel{*}{h}(X,PY),N)=\lambda g(X,PY),\quad\forall X,Y\in\Gamma(TM_{|\mathcal{U}}).
\end{equation}
\end{definition}
Note that for the case of $1$-lightlike submanifolds and lightlike hypersurfaces, $A_{N}$ is $\Gamma(S(TM))$-valued. Therefore, using (\ref{pr14}), the relation (\ref{eq4.3}) is equivalent to
\begin{equation}\label{eq4.3.1}
 A_{N}X=\lambda PX,\quad\forall X\in\Gamma(TM_{|\mathcal{U}}).
\end{equation} 
We note that by Theorem 2.5 in \cite{Dug}, page 161, $S(TM)$ is integrable. In case $\lambda=0$ on $\mathcal{U}$, we say that $S(TM)$ is totally geodesic.\vspace{0.2cm}\\
In the following, In an ambient space form, we give characterization of a family of admissible pair of screens and we prove that the Ricci tensors associated are symmetric.
%%%%%%%%%%%%%%%%%%%%%%%%%%%%%%%%
\begin{theorem}\label{th4.1}
Let $M$ be an $r$-lightlike submanifold with $r<min\{m,n\}$ or a coisotropic submanifold of a $(m+n)$-dimensional semi-Riemannian manifold of constant sectional curvature $\overline{M}(c)$ that is totally umbilical. Then, the family of admissible pair of screens reduce to the set of totally umbilical screen distributions on $M$. Also, $M$ is pointwise Osserman and for each admissible pair of screens $\{S(TM),S(TM^{\perp})\}$, the associated Ricci tensor is symmetric and $M$ is locally Einstein.
\end{theorem}
\textbf{Proof:} In this proof, we suppose that $RadTM$ is of rank $r\neq1$. For the case of $1$-lightlike submanifold, the proof is similar. For the case of lightlike hypersurface, the proof is given in \cite{Ati-Dug}, Theorem 4.3.\\
Let's consider $M$ a proper totally umbilical, that is, in (\ref{eq4.2}) we have $H^{l}\neq0$ and $H^{s}\neq0$ (in case $r$-lightlike) and $H^{l}\neq0$ (in case coisotropic). Then, using (\ref{pr17}), the induced Riemann curvature tensor takes the form
$$R(X,Y)Z=c\{g(Y,Z)X-g(X,Z)Y\}+A_{h^{l}(Y,Z)}X-A_{h^{l}(X,Z)}Y+A_{h^{s}(Y,Z)}X-A_{h^{s}(X,Z)}Y.$$
Thus, using (\ref{pr9}), (\ref{pr14}) and (\ref{eq4.2}), we obtain, for any $X,Y,Z,U\in\Gamma(TM)$,
\begin{eqnarray}\label{eq4.4}
&R(X,Y,Z,U)=c\{g(Y,Z)g(X,U)-g(X,Z)g(Y,U)\}+g(Y,Z)\overline{g}(\stackrel{*}{h}(X,PU),H^{l})\crcr
&\quad -g(X,Z)\overline{g}(\stackrel{*}{h}(Y,PU),H^{l})+\overline{g}(H^{s},H^{s})\{g(Y,Z)g(X,U) -g(X,Z)g(Y,U)\}.
\end{eqnarray}
Now let's consider an admissible pair of screens $\{S(TM),S(TM^{\perp})\}$ and let $R$ denote the associate induced Riemann curvature tensor. Since 
$R(X,Y,Z,U)=R(Z,U,X,Y)$, from (\ref{eq4.4}), we obtain, for any $X,Y,Z,U\in\Gamma(TM)$,
\begin{eqnarray}\label{eq4.5}
&&g(Y,Z)\overline{g}(\stackrel{*}{h}(X,PU),H^{l})-g(X,Z)\overline{g}(\stackrel{*}{h}(Y,PU),H^{l})\crcr
&&\qquad -g(U,X)\overline{g}(\stackrel{*}{h}(Z,PY),H^{l})+g(Z,X)\overline{g}(\stackrel{*}{h}(U,PY),H^{l})=0. 
\end{eqnarray}
Since $H^{l}\neq0$, for any $X,U\in\Gamma(TM)$, choose in (\ref{eq4.5}) the vector fields $Y$ and $Z$ such that $g(Y,Z)=1$ and $g(X,Z)=0$. Thus, we get
$$\overline{g}(\stackrel{*}{h}(X,PU),H^{l})=\lambda g(X,PU),\quad\forall X,U\in\Gamma(TM),$$
where $\lambda=\overline{g}(\stackrel{*}{h}(Z,PY),H^{l})$, that is the screen distribution $S(TM)$ is totally umbilical. \\
Conversely, suppose that for any $N\in\Gamma(ltr(TM))$, there exists a smooth function $\lambda$ on $\mathcal{U}\subset M$ such that $\overline{g}(\stackrel{*}{h}(X,PY),N)=\lambda g(X,PY),~~\forall X,Y\in\Gamma(TM)$. Then, using (\ref{pr9}), (\ref{pr14}) and (\ref{eq4.2}) we have, for any $X,Y,Z,U\in\Gamma(TM)$, 
\begin{eqnarray}\label{eq4.6}
&R(X,Y,Z,U)=c\{g(Y,Z)g(X,U)-g(X,Z)g(Y,U)\}+\lambda\{g(Y,Z)g(X,PU)\crcr
&\quad -g(X,Z)g(Y,PU)\} +\overline{g}(H^{s},H^{s})\{g(Y,Z)g(X,U) -g(X,Z)g(Y,U)\}.
\end{eqnarray}
Thus, $R$ defines an algebraic curvature tensor, that is $\{S(TM),S(TM^{\perp})\}$ is admissible.\\
Now, let $\{S(TM),S(TM^{\perp})\}$ be an arbitrary admissible pair of screens on $M$. We compute the induced Ricci curvature with respect to $\{S(TM),S(TM^{\perp})\}$. Consider the quasi-orthonormal field of frames $\{E_{1}=\xi_{1},...,E_{r}=\xi_{r},E_{r+1},...,E_{m}\}$ on $TM_{|\mathcal{U}}$. Then, for any $X,Y,Z\in\Gamma(TM)$,
\begin{eqnarray*}
R(X,Y)Z&=&c\{g(Y,Z)X-g(X,Z)Y\}+\lambda\{g(Y,Z)X-g(X,Z)Y\}\\ 
&&+g(Y,Z)A_{H^{s}}X -g(X,Z)A_{H^{s}}Y.
\end{eqnarray*}
Thus, we have, for any $X,Y\in\Gamma(TM)$,
\begin{eqnarray*}
Ric(X,Y)&=&\sum_{i=r+1}^{m}\varepsilon_{i}g(R(X,E_{i})Y,E_{i})+\sum_{i=1}^{r}g(R(X,\xi_{i})Y,N_{i})\\
&=&(c+\lambda+\overline{g}(H^{s},H^{s}))\{g(X,Y)-(m-r)g(X,Y)\}-r(c+\lambda)g(X,Y)\\
&& -g(X,Y)\sum_{i=1}^{r}\overline{g}(D^{s}(\xi_{i},N_{i}),H^{s})\\
&=&\big\{(1-m+r)(c+\lambda+\overline{g}(H^{s},H^{s}))-r(c+\lambda)
-\sum_{i=1}^{r}\overline{g}(D^{s}(\xi_{i},N_{i}),H^{s})\big\}g(X,Y).
\end{eqnarray*}
Therefore, the induced Ricci curvature is symmetric. Moreover $M$ is locally Einstein.\\
Finally, let's consider $x\in S_{u}^{+}(M)$, $(or~~x\in S_{u}^{-}(M))$, $u\in M$, $y,v\in x^{\perp}$. Then, by using (\ref{eq4.6}), we have
\begin{eqnarray*}
\tilde{g}(J_{R}(x)y,v)&=&R(y,x,x,v)\\
&=&(c+\lambda+\overline{g}(H^{s},H^{s}))g(x,x)g(y,v).
\end{eqnarray*}
Thus, the pseudo-Jacobi operator $J_{R}(x)$ and its characteristic polynomial $P_{t}(J_{R}(x))$ are given by
$$J_{R}(x)y=\pm(c+\lambda+\overline{g}(H^{s},H^{s}))Py,\quad P_{t}(J_{R}(x))=-t\{\pm(c+\lambda+\overline{g}(H^{s},H^{s}))-t\}^{m},$$
for arbitrary given admissible pair of screens. Therefore, $M$ is spacelike(timelike) pointwise Osserman, which completes the proof.
$\quad\square$
%%%%%%%%%%%%%%%%%%%%%%%%%%%%%%%%%%
\begin{proposition}\label{th4.2}
Let $M$ be an $r$-lightlike submanifold with $r<min\{m,n\}$ or a coisotropic submanifold of a $(m+n)$-dimensional semi-Riemannian manifold $(\overline{M},\overline{g})$, with induced algebraic Riemannian curvature tensor $R$. For any $x\in S_{u}^{+}(M)$, $(x\in S_{u}^{-}(M))$, $u\in M$ we have
$$trace J_{R}(x)=\sum_{i=1}^{r}\eta_{i}(R(x,\xi_{i})x)-Ric(x,x),$$
where $\eta_{i}(\cdotp)=\overline{g}(\cdotp,N_{i})$.
\end{proposition}
\textbf{Proof:} Let's consider the quasi-orthonormal basis  $\{E_{1}=\xi_{1},...,E_{r}=\xi_{r},E_{r+1},...,E_{m}\}$ of $T_{u}M$ such that $S(TM)=Span\{E_{r+1},...,E_{m}\}$. We have
\begin{eqnarray*}
trace J_{R}(x)&=&\sum_{i=1}^{m}\tilde{g}^{ii}\tilde{g}(J_{R}(x)E_{i},E_{i})\\
&=&\sum_{i=r+1}^{m}\tilde{g}^{ii}\tilde{g}(J_{R}(x)E_{i},E_{i})+
    \sum_{i=1}^{r}\tilde{g}(J_{R}(x)\xi_{i},\xi_{i})\\
&=&-\sum_{i=r+1}^{m}\tilde{g}^{ii}\tilde{g}(R(x,E_{i})x,E_{i})\\
&=&-\sum_{i=r+1}^{m}\tilde{g}^{ii}\tilde{g}(R(x,E_{i})x,E_{i})
   -\sum_{i=1}^{r}\tilde{g}(R(x,\xi_{i})x,\xi_{i})+\sum_{i=1}^{r}\tilde{g}(R(x,\xi_{i})x,\xi_{i})\\
&=&\sum_{i=1}^{r}\eta_{i}(R(x,\xi_{i})x)-Ric(x,x).\quad\square
\end{eqnarray*}

It is known (see \cite{Ga-Ku-Va}) that in semi-Riemannian case, being Osserman at a point simplifies the geometry at that point as the manifold is Einstein at that point. Moreover, if the latter is connected and of at least dimension 3, by Schur lemma (\cite{Be}), it is Einstein. For the degenerate case we have the following result.
%%%%%%%%%%%%%%%%%%%%%%%%%%%%%%%%
\begin{theorem}\label{th4.3}
Let $M$ be an $r$-lightlike submanifold with $r<min\{m,n\}$ or a coisotropic submanifold of a $(m+n)$-dimensional semi-Riemannian manifold that is Osserman at $u\in M$. If for a given admissible pair of screens $\{S(TM),S(TM^{\perp})\}$ , the associate induced curvature tensor $R$ verifies,
 $$(1)~~ R(x,\xi_{i})\xi_{j}=0,\quad R(\xi_{i},\xi_{j})x=0,~\forall x\in T_{u}M,~~1\leqslant i,j\leqslant r\qquad\qquad$$  
 $$(2)~~|\eta_{i}(R(x,\xi_{i})x)|<\mu\in\mathbb{R},\quad \forall x\in S_{u}^{+}(M)( or~~x\in S_{u}^{-}(M)),~~ 1\leqslant i\leqslant r,$$ 
then $(M,g,S(TM),S(TM^{\perp}))$ is Einstein at $u\in M$.
\end{theorem}
\textbf{Proof:} Let's consider the quasi-orthonormal basis  $\{E_{1}=\xi_{1},...,E_{r}=\xi_{r},E_{r+1},...,E_{m}\}$ of $T_{u}M$ such that $S(TM)=Span\{E_{r+1},...,E_{m}\}$. Denote by $R'$ and $g'$ the restriction on $S(TM)$ of the induced algebraic curvature tensor $R$ and the induced metric $g$ on $M$, respectively. Since the characteristic polynomial $J_{R}(x)$ is the same for any
$x\in S_{u}^{+}(M)$ $(x\in S_{u}^{-}(M))$, we have $traceJ_{R}(x)$ is bounded on $S_{u}^{+}(M)$, $S_{u}^{-}(M))$. By using Proposition \ref{th4.2}, we have for any $x\in S_{u}^{+}(M)$ $(x\in S_{u}^{-}(M))$,
$$|Ric(x,x)|\leqslant|traceJ_{R}(x)|+\sum_{i=1}^{r}|\eta_{i}(R(x,\xi_{i})x)|$$
It follows that there exist $\alpha\in\mathbb{R}$ such that $|Ric(x,x)|\leqslant\alpha~$ for any $x\in S_{u}^{+}(M)$, $(x\in S_{u}^{-}(M))$.
In particular, we have
$$|Ric'(x,x)|\leqslant\alpha,\quad\forall x\in S^{+}(S(T_{u}M))( x\in  S^{-}(S(T_{u}M) ).$$
Thus, since $(S(T_{u}M),g')$ is non-degenerate, it follows from a well known algebraic result (see \cite{Daj-Nomi}) that there exist $\lambda\in\mathbb{R}$ such that
\begin{equation}\label{eq4.7}
Ric'(x,y)=\lambda g'(x,y),\quad\forall x,y\in S(T_{u}M).
\end{equation}
Also, since $R$ is algebraic, by hypothesis we have, for any $x\in T_{u}M$, $1\leqslant i\leqslant r$,
$$Ric(\xi_{i},x)=\sum_{j=1}^{r}R(\xi_{i},\xi_{j},x,N_{j})+\sum_{j=r+1}^{m}\varepsilon_{j}R(\xi_{i},E_{j},x,E_{j})=\sum_{j=r+1}^{m}\varepsilon_{j}R( E_{j},x,E_{j},\xi_{i})=0,$$
and
$$Ric(x,\xi_{i})=\sum_{j=1}^{r}R(x,\xi_{j},\xi_{i},N_{j})+\sum_{j=r+1}^{m}\varepsilon_{j}R(x,E_{j},\xi_{i},E_{j})=-\sum_{j=r+1}^{m}\varepsilon_{j}R(x, E_{j},E_{j},\xi_{i})=0.$$
Therefore, since $g(\xi_{i},x)=0,~\forall x\in T_{u}M$, the latter and (\ref{eq4.7}) lead to
$$Ric(x,y)=\lambda g(x,y),\quad\forall x,y\in T_{u}M,$$
that is $(M,g,S(TM),S(TM^{\perp}))$ is Einstein at $u\in M$.$\quad\square$
%%%%%%%%%%%%%%%%%%%%%%%%%%%
\begin{corollary}\label{th4.4}
Let $(M,g)$ be a lightlike hypersurface of a semi-Riemannian manifold of constant sectional curvature $(\overline{M}(c),\overline{g})$. If $(M,g)$ is Osserman at $u\in M$ then for an admissible screen distribution $S(TM)$, $(M,g,S(TM))$ is Einstein at $u\in M$.
\end{corollary}
\textbf{Proof:} The induced curvature tensor $R$ on $M$ is given by
\begin{equation}\label{eq4.8}
R(X,Y)Z=c\{g(Y,Z)X-g(X,Z)Y\}+B(Y,Z)A_{N}X-B(X,Z)A_{N}Y.
\end{equation}
Since $R(x,\xi)\xi=0=R(\xi,\xi)x,~~\forall x\in T_{u}M$ and $|\eta(R(x,\xi)x)|=|c|,~~\forall x\in S_{u}^{\pm}(M)$, the result is immediate consequence of Theorem \ref{th4.3}.$\quad\square$
%%%%%%%%%%%%%%%%%%%%%%%%%%%%%%%%%%%%%%%%%%%%%%%%%%%%%%%%%%%%%%%%%%%%%%%%%%%%%%
\section{Symmetry properties on lightlike hypersurfaces}
Let's consider a lightlike hypersurface $(M,g)$ of an $(m+2)$-dimensional semi-Riemannian manifold of constant sectional curvature $(\overline{M}(c),g)$. By direct calculation, using (\ref{eq4.8}), for a normalizing pair $\{\xi,N\}$, the Ricci tensor on $M$ is given by
\begin{equation}\label{eq5.1}
Ric(X,Y)=mcg(X,Y)+B(X,Y)trA_{N}-B(A_{N}X,Y),\quad\forall X,Y\in\Gamma(TM).
\end{equation}
Now, suppose that $(M,g)$ is pointwise Osserman. In virtue of Corollary \ref{th4.4}, there exist a smooth function $\lambda$ on $\mathcal{U}\subset M$ such that
\begin{equation}\label{eq5.2}
 Ric(X,Y)=\lambda g(X,Y),\quad\forall X,Y\in\Gamma(TM_{|\mathcal{U}}).
\end{equation}
Thus, from relations (\ref{eq5.1}) and (\ref{eq5.2}), we get
\begin{equation}\label{eq5.3}
B(A_{N}\xi,X)=0,\quad\forall X\in\Gamma(TM_{|\mathcal{U}}).
\end{equation}

A lightlike hypersurface $(M,g,S(TM))$ of a semi Riemannian manifold $(\overline{M},\overline{g})$ is said to be locally symmetric \cite{Gun-Sahin-Kil}, if and only if for any $X,Y,Z,T,W\in\Gamma(TM)$ and $N\in\Gamma(tr(TM))$ the following hold
\begin{equation}\label{eq5.4} 
 g((\nabla_{W}R)(X,Y)Z,PT)=0 \quad  and \quad \overline{g}((\nabla_{W}R)(X,Y)Z,N)=0.
\end{equation}
That is $(\nabla_{W}R)(X,Y)Z=0$.
%%%%%%%%%%%%%%%%%%%%%%%%%%%%%%%%%%%
\begin{proposition}\label{th5.1}
Let $(M,g)$ be a lightlike hypersurface $(M,g)$ of a semi-Riemannian manifold of constant sectional curvature $(\overline{M}(c),g)$. The lightlike hypersurface $(M,g,S(TM))$ is locally symmetric if and only if it is totally geodesic.
\end{proposition}
\textbf{Proof:} Using relation (\ref{eq4.8}), we obtain, for any $V,X,Y\in\Gamma(TM)$,
\begin{eqnarray}\label{eq5.5}
 \overline{g}((\nabla_{V}R)(X,Y)\xi,N)&=&\overline{g}(R(X,Y)\stackrel{*}{A}_{\xi}V,N)\crcr
&=&c\{B(V,Y)\eta(X)-B(V,X)\eta(Y)\}.
\end{eqnarray}
If $(M,g,S(TM))$ is locally symmetric, by taking $Y=\xi$ into (\ref{eq5.5}), we obtain $B(V,X)=0,~~\forall V,X\in\Gamma(TM)$, that is $M$ is totally geodesic.
Conversely if $M$ is totally geodesic, since $R=\overline{R}_{|TM}$, we obtain  $\nabla_{V}R(X,Y)Z=\overline{\nabla}_{V}\overline{R}(X,Y)Z=0,~~\forall V,X,Y,Z\in\Gamma(TM)$. This complete the proof.$\quad\square$

In what follows, we consider curvature operator on a smooth manifold defined by
\begin{equation}\label{eq5.6} 
R(X,Y)=\nabla_{X}\nabla_{Y}-\nabla_{Y}\nabla_{X}-\nabla_{[X,Y]}.
\end{equation}
A lightlike hypersurface $(M,g,S(TM))$ of a semi-Riemannian manifold $(\overline{M},\overline{g})$ is said to be semi-symmetric if the following condition is stisfied (see \cite{Sahin})  
\begin{equation}\label{eq5.7} 
 (R(V_{1},V_{2})\cdot R)(X,Y,Z,T)=0\quad\forall V_{1},V_{2},X,Y,Z,T\in\Gamma(TM)
\end{equation}
where $R$ is the induced Riemann curvature on $M$.
This is equivalent to
$$-R(R(V_{1},V_{2})X,Y,Z,T)- ... -R(X,Y,Z,R(V_{1},V_{2})T)=0.$$
In general the condition (\ref{eq5.7}) is not equivalent to $(R(V_{1},V_{2})\cdot R)(X,Y)Z=0$ as in the non-degenerate setting. Indeed, by direct calculation we have for any $V_{1},V_{2},X,Y,Z,T\in\Gamma(TM)$,
\begin{eqnarray}\label{eq5.8}
&&(R(V_{1},V_{2})\cdot R)(X,Y,Z,T)=\crcr
&&\qquad\qquad g((R(V_{1},V_{2})\cdot R)(X,Y)Z,T)+(R(V_{1},V_{2}).g)(R(X,Y)Z,T).
\end{eqnarray}
%%%%%%%%%%%%%%%%%%%%%%%%%%%%%
\begin{theorem}\label{th5.2}
Let $(M,g)$ be a (proper) totally umbilical lightlike hypersurface of a semi-Riemannian manifold of constant sectional curvature $(\overline{M}(c),\overline{g})$. Then for an admissible screen distribution $S(TM)$, the lightlike hypersurface $(M,g,S(TM))$ is semi-symmetric.
\end{theorem}
\textbf{Proof:} Suppose that $S(TM)$ is an admissible screen distribution of a totally umbilical lightlike hypersurface $(M,g)$. In virtue of Theorem \ref{th4.1}, the screen distribution $S(TM)$ is totally umbilical. Then, by using relations (\ref{eq4.2.1}) and (\ref{eq4.3.1}), the induced curvature tensor $R$ associated to $S(TM)$ is given by
\begin{equation}\label{eq5.9}
 R(X,Y)Z=c\{g(Y,Z)X-g(X,Z)Y\}+\rho\lambda\{g(Y,Z)PX-g(X,Z)PY\},
\end{equation}
where $\rho$ and $\lambda$ are smooth functions on a coordinate neighbourhood $\mathcal{U}\subset M$.
By straightforward calculation using (\ref{eq5.9}), we obtain 
$$(R(V_{1},V_{2})\cdot R)(X,Y,Z,T)=0,\quad\forall V_{1},V_{2},X,Y,Z,T\in\Gamma(TM),$$
that is $(M,g,S(TM))$ is semi-symmetric.$\quad\square$
%%%%%%%%%%%%%%%%%%%%%%%%%%%%%
\begin{theorem}\label{th5.3}
Let $(M,g)$ be a pointwise Osserman lightlike hypersurface of a semi-Riemannian manifold of constant sectional curvature $(\overline{M}(c),\overline{g})$ such that $A_{N}\xi$ is a non-null vector field. Then for an admissible screen distribution $S(TM)$, the lightlike hypersurface $(M,g,S(TM))$ is semi-symmetric if and only if it is totally geodesic.
\end{theorem}
\textbf{Proof:} The induced curvature satisfies relation (\ref{eq4.8}). By straightforward calculation we obtain, for any $V_{1},V_{2},X,Y,Z\in\Gamma(TM)$
\begin{eqnarray}\label{eq5.10} 
&&(R(V_{1},V_{2})\cdot R)(\xi,X,Y,Z)=c\big\{B(V_{2},Y)g(V_{1},X)g(A_{N}\xi,Z)\crcr
&&\hspace{2cm}-B(V_{1},Y)g(V_{2},X)g(A_{N}\xi,Z)-B(X,V_{1})g(V_{2},Y)g(A_{N}\xi,Z)\crcr
&&\hspace{2cm}+B(X,V_{2})g(V_{1},Y)g(A_{N}\xi,Z)-B(X,Y)g(A_{N}\xi,V_{1})g(V_{2},Z)\crcr 
&&\hspace{2cm}+B(X,Y)g(A_{N}\xi,V_{2})g(V_{1},Z)\big\}-B(V_{2},X)B(A_{N}V_{1},Y)g(A_{N}\xi,Z)\crcr
&&\hspace{2cm}+B(V_{1},X)B(A_{N}V_{2},Y)g(A_{N}\xi,Z)-B(X,A_{N}V_{1})B(V_{2},Y)g(A_{N}\xi,Z)\crcr
&&\hspace{2cm}+B(X,A_{N}V_{2})B(V_{1},Y)g(A_{N}\xi,Z)-B(X,Y)B(V_{2},Z)g(A_{N}\xi,A_{N}V_{1})\crcr
&&\hspace{2cm}+B(X,Y)B(V_{1},Z)g(A_{N}\xi,A_{N}V_{2}).
\end{eqnarray}  
So, suppose that $M$ is semi-symmetric, by taking $V_{1}=\xi$ into (\ref{eq5.10}), we get
\begin{eqnarray}
&&B(V_{2},X)B(A_{N}\xi,Y)g(A_{N}\xi,Z)+B(X,A_{N}\xi)B(V_{2},Y)g(A_{N}\xi,Z)\crcr
&&\hspace{4cm}+B(X,Y)B(V_{2},Z)g(A_{N}\xi,A_{N}\xi)=0.
\end{eqnarray}
By assymption and using Corollary \ref{th4.4}, we infer that the relation (\ref{eq5.3}) is satisfied. Thus, we obtain
\begin{equation}\label{eq5.11}
B(X,Y)B(V_{2},Z)g(A_{N}\xi,A_{N}\xi)=0,\quad\forall V_{2},X,Y,Z\in\Gamma(TM).
\end{equation}
Since $A_{N}\xi$ is non-null, by taking $V_{2}=X$ and $Z=Y$ into (\ref{eq5.11}), we infer that $B(X,Y)=0,\forall X,Y\in\Gamma(TM)$, that is $M$ is totally geodesic. Conversely if $M$ is totally geodesic, since $R=\overline{R}_{|TM}$, we obtain  $(R(V_{1},V_{2})\cdot R)(X,Y,Z,T)=(\overline{R}(V_{1},V_{2})\cdot \overline{R})(X,Y,Z,T)=0$, for any  $V_{1},V_{2},X,Y,Z,T\in\Gamma(TM)$. This complete the proof.$\quad\square$\vspace{0.2cm}\\
From Proposition \ref{th5.1} and Theorem \ref{th5.3}, we have the following result.
\begin{corollary}\label{th5.4}
Let $(M,g)$ be a pointwise Osserman lightlike hypersurface of a semi-Riemannian manifold of constant sectional curvature $(\overline{M}(c),\overline{g})$ such that $A_{N}\xi$ is a non-null vector field. Then for an admissible screen distribution $S(TM)$, the lightlike hypersurface $(M,g,S(TM))$ is locally symmetric if and only if it is semi-symmetric.
\end{corollary}
%%%%%%%%%%%%%%%%%%%%%%%%%%%%%%%%%%%%%%
A lightlike submanifold $M$ of a semi-Riemannian manifold $\overline{M}$ is said to be Ricci semi-symmetric if the following condition is satisfied
\begin{equation}\label{eq5.12}
R(V_{1},V_{2}\cdot Ric)(X,Y)=0,~~\forall V_{1},V_{2},X,Y\in\Gamma(TM),
\end{equation}
where $R$ and $Ric$ are induced Riemann curvatrure and Ricci tensor on $M$, respectively. The latter condition is equivalent to
$$-Ric(R(V_{1},V_{2})X,Y)-Ric(X,R(V_{1},V_{2})Y)=0$$
%%%%%%%%%%%%%%%%%%%%%%%%%%%%%%%%%%%%
\begin{theorem}\label{th5.5}
Let $(M,g)$ be a pointwise Osserman lightlike hypersurface of a semi-Riemannian manifold of constant sectional curvature $(\overline{M}(c),\overline{g})$. Then for an admissible screen distribution $S(TM)$, the lightlike hypersurface $(M,g,S(TM))$ is Ricci semi-symmetric.
\end{theorem}
\textbf{Proof:} Let $R$ be an induced curvature tensor associated to the admissible screen distribution $S(TM)$.  By virtue of Corollary \ref{th4.4}, there exists a smooth function $\lambda$ on $\mathcal{U}\subset M$ such that $Ric(X,Y)=\lambda g(X,Y),~~\forall X,Y\in\Gamma(TM_{|\mathcal{U}})$. Thus, since $R$ is an algebraic curvature tensor, we have, for any $V_{1},V_{2},X,Y\in\Gamma(TM)$
\begin{eqnarray*}
(R(V_{1},V_{2})\cdot Ric)(X,Y)&=&-Ric(R(V_{1},V_{2})X,Y)-Ric(X,R(V_{1},V_{2})Y)\\
 &=&-\lambda R(V_{1},V_{2},X,Y)-\lambda R(V_{1},V_{2},Y,X)\\
&=&0,
\end{eqnarray*}
that is $(M,g,S(TM))$ is Ricci semi-symmetric.$\quad\square$

%%%%%%%%%%%%%%%%%%%%%%%%%%%%%%%%%%%%%%%%%%%%%%%%%%%%%%%%%%%%%%%%%%%%%%%%%%%%%%%

\end{document}